\documentclass[a4paper,10pt]{article}
\usepackage{amsthm}
\usepackage{amsmath}
\usepackage{amsfonts}
\usepackage{amssymb}

\theoremstyle{plain}\newtheorem{theorem}{Theorem}[section]
\theoremstyle{definition}\newtheorem{numero}[theorem]{}

\def\ov #1{\overline{#1}}
\def\BbbP{{\mathbb P}}
\def\BbbK{{\mathbb K}}
\def\BbbT{{\mathbb T}}
\def\dd{{\lower4pt\hbox{$\mskip-2mu\shortmid\mskip-2mu$}}}

\author{Maurizio Cailotto \thanks{\small Dip.Matem., Padova University, Italy - \tt maurizio@math.unipd.it}} 
\date{}
\begin{document}
\title{Classification of compact real surfaces: \\ a quick proof. }
\maketitle

\begin{abstract}
The aim of this note is to give a quick algebraic proof 
of (the combinatorial part of) the classification theorem for compact real surfaces, 
whose classical proofs (as in the Massey book and 
in the Conway ZIP proof) are based on surgery (and pictures) 
and look more intricate. 
\end{abstract}
\bigskip


\section{Algebraic rules} 
\numero
Let consider a word where all letters compare exactly twice 
(eventually inverted): sometime called a combinatorial polygon, 
it is the start point for the classification theorem of (connected) compact real surfaces. 
We will use roman letters $a,b,c,\dots$ for the edges, 
and greek letters $\alpha,\beta,\gamma,\dots$ for sequences of letters. 

\numero
Let define the following algebraic rules: 
\begin{description}
\item{$(1)$}
$\alpha a\sim a\alpha$ (cyclicity); 
\item{$(2)$}
$\alpha \sim \ov\alpha$ (inversion: $\ov{\alpha a}=\ov a\;\ov\alpha$); 
\item{$(3)$}
$\alpha a\ov a\sim\alpha$ (cancelation: the empty word represents the sphere); 
\item{$(4)$}
$\alpha a\beta\beta' \ov a\sim \alpha a\beta'\beta  \ov a$
(rule for discord letters); 
\item{$(5)$}
$\alpha a\beta a\sim \alpha\ov\beta aa$
(rule for concord letters); 
\end{description}

\numero
The first three rules are clear, 
the last two are the elementary surgery: 
cut along $a\beta$ and past along $a$ 
(and rename $a$ the new edge). 
All rules can be interpreted also in terms of relation in a free group: 
the quotients of a free group generated by the letters subject to 
equivalent relations (words) are isomorphic (use $c=a\beta$, 
then rename $a$ for $c$ in the new relation). 

\numero
Let notice the following immediate consequences: 
\begin{description}
\item{$(5')$}
using inversion, $(5)$ gives also 
$\alpha a\beta a\sim \alpha aa\ov\beta$, 
and, using cyclicity, generalizes to 
$
\alpha a\beta\beta' a\sim 
\alpha a\beta' a\ov\beta\sim 
\alpha \ov{\beta'} a\beta a\sim 
\alpha\ov{\beta'} aa\ov\beta. 
$ 
\item{$(6)$} 
using cyclicity both rules $(4)$ and $(5')$ can be written 
with invariant tails: 
$\alpha a\beta\beta' \ov a\gamma\sim 
\alpha a\beta'\beta  \ov a\gamma$ 
and 
$
\alpha a\beta\beta' a\gamma\sim 
\alpha\ov{\beta'} aa\ov\beta\gamma 
$. 
\item{$(7)$} 
blocks as $aa$ and $ab\ov a\ov b$ can be freely moved: 
$\alpha aa\beta\gamma\sim \alpha\beta aa\gamma$ 
(passing though $\alpha a\ov\beta a\gamma$), 
and 
$\alpha ab\ov a\ov b\beta\gamma\sim \alpha\beta ab\ov a\ov b\gamma$ 
(by 
\def\dd{{\lower4pt\hbox{$\mskip-2mu\shortmid\mskip-2mu$}}}
$$ 
\alpha\d{} ab\ov a\d{}\ov b\beta\dd{} \sim  
\beta ab\d{}\ov a\dd{}\alpha \d{}\ov b \sim  
\beta a\d{}b\dd{}\alpha\d{}\ov a \ov b \sim 
\beta a\dd{}\alpha \d{}b\ov a \ov b\d{} \sim  
\alpha\beta ab\ov a \ov b   
$$ 
then use cyclicity; 
we indicate with the ``underdots" the discord letters to which the rule is applied). 
\end{description}

\numero
Finally, putting as usual $\BbbP=aa$, $\BbbT=ab\ov a\ov b$ and $\BbbK=ab\ov a b$, 
and considering the connected sum $\#$ as concatenation of words with disjoint letters, 
we have the following examples:  
\begin{description}
\item{$(i)$} 
$\BbbP\#\BbbP=aabb\sim ab\ov a b=\BbbK$, 
\item{$(ii)$} 
$\BbbP\#\BbbT=aabc\ov b\ov c\sim acbabc\sim acbb\ov ac\sim\BbbP\#\BbbK$, 
therefore also $\sim\BbbP\#\BbbP\#\BbbP$. 
\item{$(iii)$} 
using $(7)$ on central pairs we have 
$a_1a_2\cdots a_na_n\cdots a_2a_1\sim\BbbP^{\# n}$; 
\item{$(iv)$} 
$a_1a_2\cdots a_{n-1}a_n\ov{a_1}\,\ov{a_2}\cdots\ov{a_{n-1}}a_n\sim 
a_1a_2\cdots a_{n-1}a_{n-1}\cdots a_2a_1a_na_n\sim 
\BbbP^{\# n}$; 
\item{$(v)$} 
$
a_1\d{}a_2\cdots a_{n-1}\dd{}a_n\d{}\ov{a_1}\,\ov{a_2}\cdots\ov{a_{n-1}}\,\ov{a_n}\sim 
a_1a_n\d{}a_2\cdots a_{n-1}\dd{}\ov{a_1}\,\ov{a_2}\cdots\ov{a_{n-1}}\,\d{}\ov{a_n}\sim \\
a_1a_n\ov{a_1}\,\ov{a_2}\cdots\ov{a_{n-1}}a_2\cdots a_{n-1}\,\ov{a_n}\sim 
\BbbT\#(a_2\cdots a_{n-1}\ov{a_2}\cdots\ov{a_{n-1}})\sim 
\BbbT^{\#[n/2]}
$ 
by induction. 
\end{description}

\section{The theorem and its proof.}

\numero
The classification theorem asserts that connected compact real surfaces 
can be reduced to: the sphere, connected sums of tori, connected sums 
of real projective plans. 
Starting with a combinatorial polygon, the quick proof is the following: 
\begin{description}
\item{$(a)$} 
using $(5)$, $(6)$ and $(7)$ we can associate concord pairs of letters, 
and hive off all occurrences of projective plans; 
so it remains a word with only discord letters; 
\item{$(b)$} 
using $(4)$, $(6)$ and $(7)$, for any pair of discord letters which 
are separated by another pair, we can associate: 
$$ 
\alpha a\d{}\mskip 1mu\beta\dd{} b\gamma \d{}\ov a\delta\ov b\varepsilon \sim 
\alpha ab\d{}\gamma \beta\dd{} \ov a\delta\d{}\ov b\varepsilon \sim 
\alpha \d{}ab\ov a\d{}\delta\gamma \beta\dd{}\ov b\varepsilon \sim 
ab \ov a\ov b\varepsilon \alpha \delta\gamma \beta 
$$ 
and hive off all occurrences of tori; 
notice that this step does not change the (discord) character of the pairs; 
\item{$(c)$} 
as usual, if there are occurrences of projective plans, 
then all occurrences of tori can be changed in couples of projective plans by $(ii)$; 
\item{$(d)$} 
finally, it remains a word with only discord letters 
and such that any discord pair of letters is not separated 
by another one; 
using (descending) induction on the number of letters between a pair, 
we will arrive at a pair $a\ov a$ and use rule $(3)$ to cancel out
all such occurrences. 
\end{description}

\numero
The proof shows also that a combinatorial polygon gives an 
orientable surface if and only if there are no concord letters in (any of) its word. 

\numero
The classification theorem for connected compact real surfaces with boundary, 
that is two such surfaces are isomorphic if and only if the corresponding 
connected compact real surfaces (obtained by gluing a $2$-cell on any component 
of the boundary) are isomorphic, can be recovered with the same computations, 
assuming the new rule that adjacent single letters can be glue into one. 
In fact, starting with a polygon (word) having also single edges (letters), 
performing steps $(a)$, $(b)$, $(c)$ as before, we are reduced to discuss 
the case in which only discord pairs of letters compares, 
never separated by other pairs, with possibly single letters. 
If there are no pairs, all reduced to a single letter, which represents a hole. 
Otherwise, using induction on the number of letters between a (discord) pair 
we will arrive at the cases $a\ov a$ and $a x\ov a$ where $x$ is a single letter; 
the first case is canceled as before, the second one can be hive off: 
any one of such sequences represents a hole for the surface. 


\end{document}